\documentclass{amsart}
\usepackage{amsmath,amsfonts,latexsym,amssymb}

\newtheorem{theorem}{Theorem}
\newtheorem{lemma}{Lemma}

\theoremstyle{remark}
\newtheorem{remark}{Remark}
\newtheorem*{ack}{Acknowledgements}

\begin{document}

\markboth{Ritabrata Munshi}{Shifted convolution sums for $GL(3)\times GL(2)$}
\title[Shifted convolution sums for $GL(3)\times GL(2)$]{Shifted convolution sums for $GL(3)\times GL(2)$}

\author{Ritabrata Munshi}   
\address{School of Mathematics, Tata Institute of Fundamental Research, 1 Dr. Homi Bhabha Road, Colaba, Mumbai 400005, India.}     
\email{rmunshi@math.tifr.res.in}

\begin{abstract}
For the shifted convolution sum  
$$
D_h(X)=\sum_{m=1}^\infty\lambda_1(1,m)\lambda_2(m+h)V\left(\frac{m}{X}\right)
$$    
where $\lambda_1(1,m)$ are the Fourier coefficients of a $SL(3,\mathbb Z)$ Maass form $\pi_1$, and $\lambda_2(m)$ are those of a $SL(2,\mathbb Z)$ Maass or holomorphic form $\pi_2$, and $1\leq |h| \ll X^{1+\varepsilon}$, we establish the bound
$$
D_h(X)\ll_{\pi_1,\pi_2,\varepsilon} X^{1-\frac{1}{20}+\varepsilon}.
$$
The bound is uniform with respect to the shift $h$.
\end{abstract}

\subjclass{11F66, 11M41}
\keywords{Maass forms, shifted convolution sum, circle method}

\maketitle


\section{Introduction}
\label{intro}

The shifted convolution sum 
$$
\mathop{\sum}_{m=1}^\infty\lambda_1(m)\lambda_2(m+h)V\left(\frac{m}{X}\right), 
$$ 
with $GL(2)$ Fourier coefficients $\lambda_i(m)$, has been investigated extensively by several authors since Selberg's seminal paper \cite{Se}. Non-trivial bound of this sum often has deep implications, e.g. subconvexity and equidistribution (QUE) (see \cite{B}, \cite{DFI-1}, \cite{DFI-2}, \cite{H}, \cite{HaM}, \cite{Ho}, \cite{HM}, \cite{J-2}, \cite{KMV}, \cite{LS}, \cite{M}, \cite{S}). 

In this paper we will consider a higher rank analogue - 
$$
D_h(X):=\sum_{m=1}^\infty\lambda_1(1,m)\lambda_2(m+h)V\left(\frac{m}{X}\right)
$$    
where $\lambda_1(1,m)$ are the Fourier coefficients of a $SL(3,\mathbb Z)$ Hecke-Maass cusp form $\pi_1$, and $\lambda_2(m)$ are those of a $SL(2,\mathbb Z)$ Hecke-Maass or Hecke holomorphic cusp form $\pi_2$. We will take $V$ to be smooth and compactly supported in $[1,2]$. Also we will take $0\leq h\ll X^{1+\varepsilon}$, as for larger shifts the trivial bound most often suffices. (Pitt \cite{P} has considered a similar sum with $\tau_3(m)$ in place of the Fourier coefficients $\lambda_1(1,m)$.) Applying Cauchy and estimates coming from Rankin-Selberg theory we obtain the following (trivial) estimate
\begin{align}
\label{trivial}
D_h(X)\ll_{\pi_1,\pi_2} X^{1+\varepsilon}.
\end{align}
Our main theorem gives a nontrivial power saving over this estimate.
\begin{theorem}
For $0\leq h \ll X^{1+\varepsilon}$, we have
$$
D_h(X)\ll_{\pi_1,\pi_2,\varepsilon} X^{1-\frac{1}{20}+\varepsilon}.
$$
\end{theorem}

As in the case of the $GL(2)$ shifted convolution sum we first apply the circle method to detect the shift using additive harmonics, and then apply Voronoi summation formula. However, unlike the $GL(2)$ case, this does not solve the problem. We are left with a complicated expression (see \eqref{approx-pw-3}), involving higher dimensional Kloosterman-type character sums (see \eqref{main-char}). Assuming square-root cancellation in the character sum one can show that we are just at the threshold, and any saving in the sum of the character sums will yield a non-trivial bound. However, except in the case of the zero shift $h=0$, it is not clear how to obtain extra cancellation. We resolve this issue by adopting Jutila's variation of the circle method with an important new input - factorizable moduli (see Section~\ref{mod-set-choice} and Remark \ref{imp-remark}). This seemingly simple idea has other important applications. In \cite{Mu} we apply this idea to several subconvexity problems. 

In the theorem, the dependence of the implied constant on the conductors of $\pi_1$ and $\pi_2$ can be explicitly given, though we do not try to do this here. In fact, in some of the estimates that we use, e.g. Lemma \ref{ram-on-av}, the implied constants are independent of the conductor of the form. Moreover it is not necessary to assume that $\pi_2$ is of full level. The same bound holds for general $\pi_2$. It is expected that extra cancellation can be obtained by averaging over $h$. However we will not take up this issue in this paper.

\ack
I thank Valentin Blomer, Tim Browning, Gergely Harcos and Matthew Young for their helpful comments. 


\section{Preliminaries}
\label{prelim}

\subsection{Preliminaries on $SL(3,\mathbb Z)$ Maass forms}

We shall first recall some basic facts about $SL(3,\mathbb Z)$ automorphic forms. Our need is minimal and, in fact, the Voronoi summation formula (of Miller and Schmid \cite{MS}, and Goldfeld and Li \cite{GL}) is all that we will be using. Suppose $\pi_1$ is a Maass form of type $(\nu_1,\nu_2)$ for $SL(3,\mathbb Z)$ which is an eigenfunction of all the Hecke operators with Fourier coefficients $\lambda_1(m_1, m_2)$, normalized so that $\lambda_1(1, 1)=1$ (for details see Goldfeld's book \cite{G}). We introduce the Langlands parameters $(\alpha_1, \alpha_2, \alpha_3)$, defined by
$$
\alpha_1=-\nu_1-2\nu_2+1,\;\;\;\alpha_2=-\nu_1+\nu_2\;\;\;\text{and}\;\;\;\alpha_3=2\nu_1+\nu_2-1.
$$
The Ramanujan-Selberg conjecture predicts that $|\text{Re}(\alpha_i)|=0$, and from the work of Jacquet and Shalika we at least know that $|\text{Re}(\alpha_i)|<\frac{1}{2}$. 

Let $g$ be a compactly supported function on $(0,\infty)$, and let $\tilde g(s)=\int_0^\infty g(x)x^{s-1}dx$ be the Mellin transform. For $\sigma>-1+\max\{-\text{Re}(\alpha_1),-\text{Re}(\alpha_2),-\text{Re}(\alpha_3)\}$ and $\ell=0,1$ define
\begin{align}
\label{gl}
G_{\ell}(y)=\frac{1}{2\pi i}\int_{(\sigma)}(\pi^3 y)^{-s}\frac{\Gamma\left(\frac{1+s+\alpha_1+\ell}{2}\right)\Gamma\left(\frac{1+s+\alpha_2+\ell}{2}\right)\Gamma\left(\frac{1+s+\alpha_3+\ell}{2}\right)}{\Gamma\left(\frac{-s-\alpha_1+\ell}{2}\right)\Gamma\left(\frac{-s-\alpha_2+\ell}{2}\right)\Gamma\left(\frac{-s-\alpha_3+\ell}{2}\right)}\tilde g(-s)ds
\end{align}
and set
$$
G_+(y)=\frac{1}{2\pi^{3/2}}\left(G_{0}(y)-iG_{1}(y)\right),\;\;\;\text{and}\;\;\;G_-(y)=\frac{1}{2\pi^{3/2}}\left(G_{0}(y)+iG_{1}(y)\right).
$$
\begin{lemma}
Let $g$ be a compactly supported smooth function on $(0,\infty)$, we have
\begin{align}
\label{voronoi3}
\sum_{m=1}^\infty \lambda_1(1,m)e_q\left(am\right)g(m)=&q\sum_{m_1|q}\sum_{m_2=1}^\infty \frac{\lambda_1(m_2,m_1)}{m_1m_2}S(\bar a, m_2; q/m_1)G_+\left(\frac{m_1^2m_2}{q^3}\right)\\
\nonumber &+q\sum_{m_1|q}\sum_{m_2=1}^\infty \frac{\lambda_1(m_2,m_1)}{m_1m_2}S(\bar a, -m_2; q/m_1)G_-\left(\frac{m_1^2m_2}{q^3}\right),
\end{align}
where $(a,q)=1$, and $\bar{a}$ denotes the multiplicative inverse of $a\bmod{q}$. Also $e_q(x)=e^{2\pi ix/q}$.
\end{lemma}

\begin{remark}
If $g$ is supported in $[X,2X]$, satisfying $x^jg^{(j)}(x)\ll_j 1$, then the sums on the right hand side of \eqref{voronoi3} are essentially supported on $m_1^2m_2\ll q^3(qX)^{\varepsilon}/X$ (where the implied constant depends on the form $\pi_1$ and $\varepsilon$). The contribution from the terms with $m_1^2m_2\gg q^3(qX)^{\varepsilon}/X$ is negligibly small (i.e. $O((qX)^{-N})$ for any $N>0$). This follows by estimating the integral $G_{\ell}(y)$ by shifting the contour to the right. For smaller values of $m_1^2m_2$ we shift the contour to left upto $\sigma=-\frac{1}{2}$ (using the result of Jacquet and Shalika) to obtain
\begin{align}
\label{g-bound}
G_{\pm}\left(\frac{m_1^2m_2}{q^3}\right)\ll \sqrt{\frac{Xm_1^2m_2}{q^3}}.
\end{align}
\end{remark}

The following lemma is well-known.
\begin{lemma}
\label{ram-on-av}
We have
$$
\sum_{n\leq x}|\lambda_1(1,n)|^2\ll x^{1+\varepsilon},
$$
where the implied constant depends on the form $\pi_1$ and $\varepsilon$.
\end{lemma}

\subsection{Preliminaries on $SL(2,\mathbb Z)$ Maass forms}
 
Next we shall recall the Voronoi summation formula for $SL(2,\mathbb Z)$ automorphic forms. For the sake of exposition we only present the case of Maass forms. The case of holomorphic forms is just similar (or even simpler). Furthermore, for technical simplicity, we only restrict to the case of full level. Let $\pi_2$ be a Maass cusp form with Laplace eigenvalue $\frac{1}{4}+\mu^2\geq 0$, and with Fourier expansion
$$
\sqrt{y}\sum_{n\neq 0}\lambda_2(n)K_{i\mu}(2\pi|n|y)e(nx).
$$     
We will use the following Voronoi type summation formula (see Meurman \cite{Me}).  
\begin{lemma}
Let $h$ be compactly supported smooth function on $(0,\infty)$. We have
\begin{align}
\label{voronoi2}
\sum_{n=1}^\infty \lambda_2(n)e_q\left(an\right)h(n)=\frac{1}{q}\sum_{\pm}\sum_{n=1}^\infty \lambda_2(\mp n)e_q\left(\pm\bar{a}n\right)H^{\pm}\left(\frac{n}{q^2}\right)
\end{align}
where $\bar{a}$ is the multiplicative inverse of $a\bmod{q}$, and
\begin{align*}
H^-(y)=&-\frac{\pi}{\cosh \pi\mu}\int_0^\infty h(x)\{Y_{2i\mu}+Y_{-2i\mu}\}\left(4\pi\sqrt{xy}\right)dx\\
H^+(y)=&4\cosh \pi\mu\int_0^\infty h(x)K_{2i\mu}\left(4\pi\sqrt{xy}\right)dx.
\end{align*}
\end{lemma}

\begin{remark}
If $h$ is supported in $[Y,2Y]$, satisfying $y^jh^{(j)}(y)\ll_j 1$, then the sums on the right hand side of \eqref{voronoi2} are essentially supported on $n\ll q^2(qY)^{\varepsilon}/Y$ (where the implied constant depends on the form $\pi_2$ and $\varepsilon$). The contribution from the terms with $n\gg q^2(qY)^{\varepsilon}/Y$ is negligibly small. For smaller values of $n$ we will use the trivial bound $H^{\pm}(n/q^2)\ll Y$.
\end{remark}


\section{Applying the circle method}
\label{cm}

\subsection{A version of the circle method}
We will be using a variant of the circle method, with overlapping intervals, which has been investigated by Jutila (\cite{J-1}, \cite{J-2}). For any set $S \subset \mathbb R$, let $\mathbb I_S$ denote the associated characteristic function, i.e. $\mathbb I_S(x)=1$ for $x\in S$ and $0$ otherwise. For any collection of positive integers $\mathcal Q \subset [1,Q]$ (which we call the set of moduli), and a positive real number $\delta$ in the range $Q^{-2}\ll \delta \ll Q^{-1} $, we define the function
$$
\tilde I_{\mathcal Q,\delta} (x)=\frac{1}{2\delta L}\sum_{q\in\mathcal Q}\;\sideset{}{^\star}\sum_{a\bmod{q}}\mathbb{I}_{[\frac{a}{q}-\delta,\frac{a}{q}+\delta]}(x),
$$
where $L=\sum_{q\in\mathcal Q}\phi(q)$. This is an approximation for $\mathbb I_{[0,1]}$ in the following sense: 
\begin{lemma}
\label{jutila-lemma}
We have
\begin{align}
\label{jutila}
\int_0^1\left|1-\tilde I_{\mathcal Q,\delta}(x)\right|^2dx\ll \frac{Q^{2+\varepsilon}}{\delta L^2}.
\end{align}
\end{lemma}

To prove this, let
$$
a_n:=\int_0^1\tilde I_{\mathcal Q,\delta}(x)e(-nx)dx=\frac{1}{2\delta L}\sum_{q\in\mathcal Q}\;\left[\sideset{}{^\star}\sum_{a\bmod{q}}e_q(-an)\right]\int_{-\delta}^{\delta}e(-nx)dx
$$
be the $n$-th Fourier coefficient of $\tilde I_{\mathcal Q,\delta}(x)$. The sum over $a\bmod{q}$ is the Ramanujan sum $c_q(n)$, which can be bounded as $|c_q(n)|\leq \sum_{d|(n,q)}d$. Clearly $a_0=1$. For $n\neq 0$, by estimating the integral trivially we get that
\begin{align}
\label{bd-fc}
|a_n|\ll \frac{1}{L}\sum_{q\in\mathcal Q}\;\sum_{d|(n,q)}d\ll \frac{1}{L}\sum_{q\leq Q}\;\sum_{d|(n,q)}d\ll \frac{1}{L}\sum_{\substack{d|n\\d\leq Q}}d\sum_{\substack{q\leq Q\\q\equiv 0\bmod{d}}}1\ll \frac{Q|n|^{\varepsilon}}{L}.
\end{align}
On the other hand we can also bound the integral by $\ll |n|^{-1}$, and obtain
\begin{align}
\label{bd-fc-2}
|a_n|\ll \frac{Q|n|^{\varepsilon}}{\delta L|n|}.
\end{align}
By Parseval we have
$$
\int_0^1\left|1-\tilde I_{\mathcal Q,\delta}(x)\right|^2dx=\int_0^1\left|\sum_{n\neq 0}a_ne(nx)\right|^2dx=\sum_{n\neq 0}|a_n|^2.
$$
Now to conclude the lemma we apply \eqref{bd-fc} for $|n|\leq \delta^{-1}$, and \eqref{bd-fc-2} for $|n|> \delta^{-1}$.


\subsection{Setting up the circle method}

Let $W$ be a smooth function supported in $[1/2,3]$ satisfying $W(x)=1$ for $x\in [1,2]$, and let $Y=X+h$. (We note that $X\leq Y \ll X^{1+\varepsilon}$.) Then we have
\begin{align*}
D_h(X)=&\mathop{\sum\sum}_{m,n=1}^\infty\lambda_1(1,m)\lambda_2(n)V\left(\frac{m}{X}\right)W\left(\frac{n}{Y}\right)\delta(n,m+h)\\
=&\int_0^1e(xh)\left[\mathop{\sum}_{m=1}^\infty\lambda_1(1,m)e(xm)V\left(\frac{m}{X}\right)\right]\left[\mathop{\sum}_{n=1}^\infty\lambda_2(n)e(-xn)W\left(\frac{n}{Y}\right)\right]dx.
\end{align*}
Let $\mathcal Q$ be a collection of moduli of size $Q$. Suppose $|\mathcal Q|\gg Q^{1-\varepsilon}$, so that $L=\sum_{q\in\mathcal Q}\phi(q)\gg Q^{2-\varepsilon}$. Let
\begin{align}
\label{approx-0}
\tilde D_h(X,\mathcal Q):=\int_0^1\tilde I_{\mathcal Q,\delta}(x)e(xh)\left[\mathop{\sum}_{m=1}^\infty\lambda_1(1,m)e(xm)V\left(\frac{m}{X}\right)\right]\left[\mathop{\sum}_{n=1}^\infty\lambda_2(n)e(-xn)W\left(\frac{n}{Y}\right)\right]dx.
\end{align}
For convenience we will use the short hand notation $\tilde D_h(X)$ in place of $\tilde D_h(X,\mathcal Q)$. Using the definition of the approximating function $\tilde I_{\mathcal Q,\delta}(x)$, we get
\begin{align}
\label{approx}
\tilde D_h(X)=\frac{1}{L}\sum_{q\in\mathcal Q}\;\sideset{}{^\star}\sum_{a\bmod{q}}e_q(ah)\mathop{\sum\sum}_{m,n\in \mathbb Z}\lambda_1(1,m)\lambda_2(n)e_q(a(m-n))F(m,n)
\end{align}
where
$$
F(x,y)=V\left(\frac{x}{X}\right)W\left(\frac{y}{Y}\right)\frac{1}{2\delta}\int_{-\delta}^{\delta}e(\alpha(x-y))d\alpha.
$$
We choose $\delta=Y^{-1}$ so that 
$$
\frac{\partial^{i+j}}{\partial^i x\partial^j y}F(x,y)\ll_{i,j} \frac{1}{X^iY^j}.
$$
In circle method we approximate $D_h(X)$ by $\tilde D_h(X)$, and then try to estimate the latter sum. Lemma~\ref{jutila-lemma} gives a way to estimate the error of replacing $D_h(X)$ by $\tilde D_h(X)$. More precisely we have
\begin{align}
\label{error}
\left|D_h(X)-\tilde D_h(X)\right|\ll \int_0^1\left|\mathop{\sum}_{m=1}^\infty\lambda_1(1,m)e(xm)V\left(\tfrac{m}{X}\right)\right|\left|\mathop{\sum}_{n=1}^\infty\lambda_2(n)e(-xn)W\left(\tfrac{n}{Y}\right)\right|\left|1-\tilde I(x)\right|dx,
\end{align}

Using the well-known point-wise uniform bound 
$$
\mathop{\sum}_{n=1}^\infty\lambda_2(n)e(-xn)W\left(\frac{n}{Y}\right)\ll_{\pi_2} Y^{\frac{1}{2}+\varepsilon}
$$
it follows that the right hand side of \eqref{error} is bounded by
\begin{align*}
\ll Y^{\frac{1}{2}+\varepsilon}\int_0^1\left|\mathop{\sum}_{m=1}^\infty\lambda_1(1,m)e(xm)V\left(\tfrac{m}{X}\right)\right|\left|1-\tilde I(x)\right|dx.
\end{align*}
Now we apply Cauchy and Lemma \ref{jutila-lemma} to conclude
$$
D_h(X)=\tilde D_h(X)+O\left(\frac{\sqrt{XY}Q(XYQ)^{\varepsilon}}{\sqrt{\delta}L}\right).
$$
As the moduli set $\mathcal Q$ is such that $L\gg Q^{2-\varepsilon}$, and $\delta=Y^{-1}$, it follows that the above error term is smaller than the trivial bound \eqref{trivial}, if we choose $Q=Y^{\frac{1}{2}+\delta}$ for any $\delta>0$. Indeed with this choice we have the following:
\begin{lemma}
\label{first-lemma}
We have
\begin{align}
\label{exp1}
D_h(X)=\tilde D_h(X)+O\left(X^{1-\delta+\varepsilon}\right).
\end{align}
\end{lemma}


\section{Estimation of $\tilde D_h(X)$}
\label{gl1_twists}

\subsection{Applying the Voronoi summation formulas}
Fix $\alpha\in [-\delta,\delta]$ and set
\begin{align*}
\tilde D_{h,\alpha}(X)=\frac{1}{L}\sum_{q\in\mathcal Q}\;\sideset{}{^\star}\sum_{a\bmod{q}}e_q(ah)\mathop{\sum\sum}_{m,n\in \mathbb Z}\lambda_1(1,m)\lambda_2(n)e_q(a(m-n))V\left(\frac{m}{X}\right)W\left(\frac{n}{Y}\right)e(\alpha(m-n)),
\end{align*}
so that $\tilde D_{h}(X)=\frac{1}{2\delta}\int_{-\delta}^{\delta}\tilde D_{h,\alpha}(X)d\alpha$. Now we apply Voronoi summations on the sums over $m$ and $n$. This process gives rise to several terms as noted in Section \ref{prelim} - Lemma \ref{voronoi3} and Lemma \ref{voronoi2}. As far as our analysis is concerned all the terms are of equal complexity, and so we just focus our attention on one such term -
\begin{align}
\label{approx-pw-1}
\tilde D_{h,\alpha,1}(X)=\frac{1}{L}\sum_{q\in\mathcal Q}&\sum_{m_1|q}\sum_{m_2=1}^\infty \frac{\lambda_1(m_2,m_1)}{m_1m_2}\\
\nonumber &\times\sum_{n=1}^\infty \lambda_2(-n)\mathcal S(m_1,m_2,n,h;q)G_+\left(\frac{m_1^2m_2}{q^3}\right)H^{+}\left(\frac{n}{q^2}\right),
\end{align}
where the character sum is given by
\begin{align}
\label{main-char}
\mathcal S(m_1,m_2,n,h;q):=\sideset{}{^\star}\sum_{a\bmod{q}}e_q(ah)e_q(-\bar{a}n)S(\bar a, m_2; q/m_1). 
\end{align}
Also here we are taking
$$
g(x)=V\left(\frac{x}{X}\right)e(\alpha x),\;\;\;\text{and}\;\;\;h(y)=W\left(\frac{y}{Y}\right)e(-\alpha y).
$$
The functions $G_+$ and $H^+$ are defined in Lemma \ref{voronoi3} and Lemma \ref{voronoi2} respectively.

\begin{remark}
Suppose we establish square-root cancellation in the character sum $\mathcal S(m_1,m_2,n,h;q)$ in \eqref{approx-pw-1}. Then estimating the remaining sums trivially using the decay in the functions $G_+$ and $H^+$, we get that 
$$
\tilde D_{h,\alpha,1}(X)\ll Y^{1+2\delta} (XY)^{\varepsilon} \asymp X^{1+2\delta+\varepsilon}.
$$
(Recall that $Y\ll X^{1+\varepsilon}$.) This yields the bound $D_h(X)\ll X^{1+2\delta+\varepsilon}$, which is worse than the trivial bound by an arbitrary small power $X^{2\delta}$. 
\end{remark}

Our job now is to get a nontrivial estimate for \eqref{approx-pw-1}, beyond square-root cancellation in the character sum $\mathcal S(m_1,m_2,n,h;q)$. For $h=0$, the zero shift, the character sum $\mathcal S(m_1,m_2,n,0;q)$ can be evaluated precisely, and then one can use the large sieve inequality of Duke, Friedlander and Iwaniec \cite{DFI-3} for Kloosterman fractions to get extra cancellation on the sum over $n$ and $m$. Alternatively one can use reciprocity and then Voronoi yet again on the sum over $m_2$, to get a much better result. However for non-zero shift the character sum $\mathcal S(m_1,m_2,n,h;q)$ can not be computed explicitly, and hence it is not clear how to obtain extra cancellation. This is the main issue. We will resolve this by choosing the set of moduli $\mathcal Q$ in a specific manner to get a huge structural advantage. From now on we take $h\neq 0$.

\subsection{Choosing the moduli set $\mathcal Q$}
\label{mod-set-choice}

We choose the set of moduli $\mathcal Q$ to be the product set $\mathcal Q_1 \mathcal Q_2$, where $\mathcal Q_i$ consists of primes in the dyadic segment $[Q_i,2Q_i]$ (and not dividing $h$) for $i=1,2$, and $Q_1Q_2=Q=Y^{\frac{1}{2}+\delta}$. Also we pick $Q_1$ and $Q_2$ (whose optimal sizes will be determined later) so that the collections $\mathcal Q_1$ and $\mathcal Q_2$ are disjoint. In this case, for any given $q=q_1q_2\in \mathcal Q$, we have $m_1=1,q_1,q_2$ or $q_1q_2$ in \eqref{approx-pw-1}.

If $m_1=q$ then $\mathcal S(q,m_2,n,h;q)=S(h,-n;q)$ is the usual Kloosterman sum. The well-known Weil bound gives the square-root cancellation in this case (recall that by choice $(h,q)=1$). If $m_1=q_1$, then the character sum splits as 
$$
\mathcal S(q_1,m_2,n,h;q)=S(\bar{q_2}h,-\bar{q_2}n;q_1)\sum_{a, b\in \mathbb F_{q_2}^\times}e_{q_2}(\bar{q_1}ah-\bar{q_1}\bar{a}n+b\bar a+m_2\bar b). 
$$ 
Using Weil bound we can bound the last character sum by $q_2^{3/2}$. Square-root cancellations for such sums was established in the general case by Adolphson and Sperger \cite{AS}, Denef and Loeser \cite{DL}. To adopt their result in the context of the above special sum, we consider the Newton polyhedron $\Delta(f)$ of $f(x,y)=\bar{q_1}hx-\bar{q_1}nx^{-1}+x^{-1}y+m_2y^{-1}\in \mathbb F_{q_2}[x,y,(xy)^{-1}]$. In the generic case where $q_2\nmid nm_2$, the polyhedron $\Delta(f)$ is given by the $4$-gon in $\mathbb R^2$ with vertices $(1,0)$, $(-1,0)$, $(-1,1)$ and $(0,-1)$. Hence $\dim \Delta(f)=2$. Also it is easily seen that $f$ is nondegenerate with respect to $\Delta(f)$. Adolphson and Sperger show that this condition is sufficient to conclude that 
\begin{align}
\label{as}
\sum_{a, b\in \mathbb F_{q_2}^\times}e_{q_2}(\bar{q_1}ah-\bar{q_1}\bar{a}n+b\bar a+m_2\bar b)\ll q_2
\end{align}
in the light of the general results of Deligne. 
\begin{lemma}
\label{sq-rt-can}
For $q=q_1q_2$ a product of two primes, $m_1|q$, and $(n,q_1q_2)=1$, we have
$$
\mathcal S(m_1,m_2,n,h;q)\ll \frac{q}{\sqrt{m_1}}\sqrt{\left(\frac{q}{m_1},m_2\right)}.
$$
\end{lemma}

We will now use the above lemma to estimate the following sub-sum of \eqref{approx-pw-1}, 
\begin{align}
\label{approx-pw-2}
\tilde D_{h,\alpha,1}^{\dagger}(X)=\frac{1}{L}\sum_{q_1\in\mathcal Q_1}&\sum_{q_2\in\mathcal Q_2}\sum_{m_2=1}^\infty \frac{\lambda_1(m_2,q_1)}{q_1m_2}\\
\nonumber &\times\sum_{n=1}^\infty \lambda_2(-n)\mathcal S(q_1,m_2,n,h;q_1q_2)G_+\left(\frac{m_2}{q_1q_2^3}\right)H^{+}\left(\frac{n}{q_1^2q_2^2}\right).
\end{align}
Using Lemma \ref{sq-rt-can}, and the properties of the weight functions (see Section \ref{prelim}), it follows that upto a negligible error term (i.e. $O(X^{-N})$ for any $N>0$) we have
\begin{align*}
\tilde D_{h,\alpha,1}^{\dagger}(X)\ll \frac{X^{\frac{3}{2}}}{LQ_1\sqrt{Q_2}}\sum_{q_1\in\mathcal Q_1}&\sum_{q_2\in\mathcal Q_2}\sum_{m_2\ll \frac{Q_2}{Q_1}X^{2\delta+\varepsilon}}\frac{|\lambda_1(m_2,q_1)|}{\sqrt{m_2}}\sqrt{(q_2,m_2)}\sum_{|n|\ll X^{2\delta+\varepsilon}} |\lambda_2(-n)|.
\end{align*}
(We will take $Q_1, Q_2\gg X^{2\delta+\varepsilon}$, see the remark below, so that the coprimality condition $(q_1,n)=1$ of the lemma is satisfied.) Using Cauchy inequality and the Rankin-Selberg theory the sum over $n$ is bounded by $X^{2\delta+\varepsilon}$. To the sum over $m_2$, we apply Cauchy to get
\begin{align*}
\sum_{m_2\ll \frac{Q_2}{Q_1}X^{2\delta+\varepsilon}}\frac{|\lambda_1(m_2,q_1)|}{\sqrt{m_2}}\sqrt{(q_2,m_2)}\leq \left[\sum_{m_2\ll \frac{Q_2}{Q_1}X^{2\delta+\varepsilon}}|\lambda_1(m_2,q_1)|^2\right]^\frac{1}{2}\left[\sum_{m_2\ll \frac{Q_2}{Q_1}X^{2\delta+\varepsilon}}\frac{(q_2,m_2)}{m_2}\right]^\frac{1}{2}.
\end{align*}
The last sum is clearly bounded by $X^{\varepsilon}$. To bound the middle sum we plug in the inequality
$$
|\lambda_1(m_2,q_1)|^2\leq 2|\lambda_1(m_2,1)|^2|\lambda_1(q_1,1)|^2+2|\lambda_1(m_2/q_1,1)|^2,
$$
which is a consequence of the Hecke relation, and apply Lemma \ref{ram-on-av}. It follows that upto a negligible error term we have
\begin{align*}
\tilde D_{h,\alpha,1}^{\dagger}(X)\ll \frac{X^{\frac{3}{2}+3\delta+\varepsilon}}{LQ_1^{\frac{3}{2}}}\sum_{q_1\ll Q_1}&\sum_{q_2\ll Q_2}\left(|\lambda_1(q_1,1)|+1\right)\ll \frac{X^{\frac{3}{2}+3\delta+\varepsilon}}{QQ_1^{\frac{3}{2}}}.
\end{align*}
The last inequality follows from another application of Cauchy inequality and Lemma \ref{ram-on-av}. We summarize the outcome of the above analysis in the following:
\begin{lemma}
\label{mid-lemma}
We have
\begin{align*}
\tilde D_{h,\alpha,1}^{m_1\neq 1}(X):=\frac{1}{L}\sum_{q\in\mathcal Q}&\sum_{\substack{m_1|q\\m_1\neq 1}}\sum_{m_2=1}^\infty \frac{\lambda_1(m_2,m_1)}{m_1m_2}\\
\nonumber &\times\sum_{n=1}^\infty \lambda_2(-n)\mathcal S(m_1,m_2,n,h;q)G_+\left(\frac{m_1^2m_2}{q^3}\right)H^{+}\left(\frac{n}{q^2}\right)\ll X\frac{X^{2\delta+\varepsilon}}{\min\{Q_1,Q_2\}^{\frac{3}{2}}}.
\end{align*}
\end{lemma}

\begin{remark}
Suppose  $\min\{Q_1,Q_2\}\gg X^{2\delta+\varepsilon}$, then the bound given in the above lemma is at least as good as the bound we have for the error term in \eqref{exp1}. 
\end{remark}

Applying the bound from Lemma \ref{sq-rt-can}, we can also get a satisfactory bound for the contribution from the small values of $m_2$ even when $m_1=1$. 
\begin{lemma}
\label{last-but-one}
We have
\begin{align*}
\tilde D_{h,\alpha,1}^{m_1=1}(X,M):=\frac{1}{L}\sum_{q\in\mathcal Q}&\sum_{m_2\leq M}\frac{\lambda_1(m_2,1)}{m_2}\\
\nonumber &\times\sum_{n=1}^\infty \lambda_2(-n)\mathcal S(1,m_2,n,h;q)G_+\left(\frac{m_2}{q^3}\right)H^{+}\left(\frac{n}{q^2}\right)\ll X\frac{\sqrt{M}X^{\varepsilon}}{X^{\frac{1}{4}-\frac{\delta}{2}}}.
\end{align*}
\end{lemma}

\begin{remark}
This is at least as good as the bound we have for the error term in \eqref{exp1} if $M\leq X^{\frac{1}{2}-3\delta}$. 
\end{remark}

\section{Estimation of $\tilde D_h(X)$ : Final analysis}

\subsection{Applying Cauchy and Poisson}
It remains to analyse the size of the sum
\begin{align}
\label{approx-pw-3}
\tilde D_{h,\alpha,1}^\sharp(X,M)=\frac{1}{L}\sum_{q\in\mathcal Q}\sum_{m\sim M}\frac{\lambda_1(m,1)}{m}\sum_{n=1}^\infty \lambda_2(-n)\mathcal S(1,m,n,h;q)G_+\left(\frac{m}{q^3}\right)H^{+}\left(\frac{n}{q^2}\right),
\end{align}
where $\mathcal S(1,m,n,h;q)$ is defined in \eqref{main-char}. Here $m\sim M$ means that $m$ runs over the integers in the dyadic segment $[M,2M)$. Also by Lemma \ref{last-but-one} it is enough to take $X^{\frac{1}{2}-3\delta}<M<X^{\frac{1}{2}+3\delta+\varepsilon}$. We have already noted that the square-root cancellation in the character sum is not enough for our purpose, and we need to prove cancellation in the average. To this end we will exploit heavily the factorization of the moduli set $\mathcal Q$. We have 
\begin{align*}
\tilde D_{h,\alpha,1}^\sharp(X,M)\ll \frac{1}{LM}\sum_{q_2\in\mathcal Q_2}&\sum_{n\ll X^{2\delta+\varepsilon}} |\lambda_2(-n)|\\
&\times \mathop{\sum}_{m\sim M} |\lambda_1(m,1)|\left|\sum_{q_1\in\mathcal Q_1}\mathcal S(1,m,n,h;q)G_+\left(\frac{m}{q^3}\right)H^{+}\left(\frac{n}{q^2}\right)\right|,
\end{align*}
where $q=q_1q_2$. To get rid off the Fourier coefficients, we apply Cauchy twice, Lemma \ref{ram-on-av} and its well-known $GL(2)$ version. With this we arrive at the following
\begin{align}
\label{dh}
\tilde D_{h,\alpha,1}^\sharp(X,M)\ll \frac{\sqrt{Q_2}X^{\delta+\varepsilon}}{L\sqrt{M}}\left\{\sum_{q_2\in\mathcal Q_2}\sum_{n\ll X^{2\delta+\varepsilon}}\tilde D_{h,\alpha,1}^\sharp(X,M;n,q_2)\right\}^{\frac{1}{2}},
\end{align}
where
\begin{align*}
\tilde D_{h,\alpha,1}^\sharp(X,M;n,q_2)=\mathop{\sum}_{m\in \mathbb Z}F\left(\frac{m}{M}\right)\left|\sum_{q_1\in\mathcal Q_1}\mathcal S(1,m,n,h;q_1q_2)G_+\left(\frac{m}{q_1^3q_2^3}\right)H^{+}\left(\frac{n}{q_1^2q_2^2}\right)\right|^2.
\end{align*}
Here $F$ is non-negative smooth function on $(0,\infty)$, supported on $[1/2,3]$, and such that $F(x)=1$ for $x\in [1,2]$. 

\begin{remark}
\label{imp-remark}
It is quite natural to split the set of moduli at this point. Indeed if $Q_1=1$ then we do not have enough points of summation inside the absolute value square to get more cancellation beyond the square root cancellation in the character sum $\mathcal S(1,m,n,h;q_1q_2)$. This term shows up as the diagonal contribution. On the other hand if $Q_2=1$, or $Q_1=Q$, then when we apply Poisson summation on the sum over $m$ after opening the absolute square, the size of the modulus is $Q^2=X^{1+2\delta}$, which is too large compared to the range of summation of $m$ and so the saving from Poisson is not enough. Hence the off-diagonal term is not satisfactory.    
\end{remark}

Opening the absolute square and interchanging the order of summations we get
\begin{align*}
\tilde D_{h,\alpha,1}^\sharp(X,M;n,q_2)=&\sum_{q_1\in\mathcal Q_1}\sum_{\tilde q_1\in\mathcal Q_1}H^{+}\left(\frac{n}{q_1^2q_2^2}\right)\bar H^{+}\left(\frac{n}{\tilde q_1^2q_2^2}\right)\\
&\times \mathop{\sum}_{m\in \mathbb Z}F\left(\frac{m}{M}\right)\mathcal S(1,m,n,h;q_1q_2)\bar{\mathcal S}(1,m,n,h;\tilde q_1q_2) G_+\left(\frac{m}{q_1^3q_2^3}\right)\bar G_+\left(\frac{m}{\tilde q_1^3q_2^3}\right).
\end{align*}
Applying Poisson summation on the sum over $m$ with modulus $q_1\tilde q_1q_2$, we get 
\begin{align}
\label{poisson-last}
\frac{M}{q_2}\sum_{q_1\in\mathcal Q_1}&\sum_{\tilde q_1\in\mathcal Q_1}\frac{1}{q_1\tilde q_1}H^{+}\left(\frac{n}{q_1^2q_2^2}\right)\bar H^{+}\left(\frac{n}{\tilde q_1^2q_2^2}\right)\sum_{m\in\mathbb Z}\mathcal T(n,m,h;q_1,\tilde q_1,q_2)\mathcal I(m;q_1,\tilde q_1,q_2),
\end{align}
where the character sum is given by
$$
\mathcal T(n,m,h;q_1,\tilde q_1,q_2)=\sum_{\alpha\bmod{q_1\tilde q_1q_2}}\mathcal S(1,\alpha,n,h;q_1q_2)\bar{\mathcal S}(1,\alpha,n,h;\tilde q_1q_2)e_{q_1\tilde q_1q_2}(m\alpha),
$$
and the integral is given by
$$
\mathcal I(m;q_1,\tilde q_1,q_2)=\int_{\mathbb R} F\left(x\right)G_+\left(\frac{xM}{q_1^3q_2^3}\right)\bar G_+\left(\frac{xM}{\tilde q_1^3q_2^3}\right)e_{q_1\tilde q_1q_2}(-mMx)dx.
$$
Integrating by parts repeatedly we get that the integral is negligibly small unless $|m|\ll Q_1QX^{\varepsilon}/M$. Observe that differentiating under the integral sign in \eqref{gl}, one can show that $y^jG_{\ell}^{(j)}(y)\ll_j \sqrt{yX}$. So for $|m|\ll Q_1QX^{\varepsilon}/M$ we have the bound
$$
\mathcal I(m;q_1,\tilde q_1,q_2)\ll\frac{MX}{Q^3}.
$$
The following lemma now follows from \eqref{poisson-last}.
\begin{lemma}
\label{rem-sum}
For any $N>0$, we have the bound
\begin{align*}
\tilde D_{h,\alpha,1}^\sharp(X,M;n,q_2)\ll \frac{M^2X^3}{Q_1Q^4}\sum_{q_1\in\mathcal Q_1}&\sum_{\tilde q_1\in\mathcal Q_1}\sum_{|m|\ll \frac{Q_1Q}{M}X^{\varepsilon}}|\mathcal T(n,m,h;q_1,\tilde q_1,q_2)|+X^{-N}.
\end{align*}
\end{lemma}


\subsection{Estimating the character sums}

First consider the case where $q_1\neq \tilde q_1$. Then the character sum $\mathcal T(n,m,h;q_1,\tilde q_1,q_2)$ splits into a product of three character sums with moduli $q_1$, $\tilde q_1$ and $q_2$ respectively. The sum modulo $q_1$ is given by
$$
\mathcal T_1=\sum_{\alpha\bmod{q_1}}\;\sideset{}{^\star}\sum_{\beta\bmod{q_1}}e_{q_1}(\bar{q_2}h\beta-\bar{q_2}n\bar{\beta})S(\bar{q_2}\bar{\beta},\bar{q_2}\alpha;q_1)e_{q_1}\left(\bar{\tilde q}_1\bar{q_2}m\alpha\right).
$$
Opening the Kloosterman sum and executing the sum over $\alpha$ we arrive at an explicit expression of this character sum in terms of Kloosterman sums, namely
$$
\mathcal T_1=q_1S(\bar{q_2}h,-\bar{q_2}(n+\tilde q_1\bar{m});q_1)
$$
if $(m,q_1)=1$, and $\mathcal T_1=0$ otherwise. The sum modulo $\tilde q_1$, which we denote by $\tilde{\mathcal T}_1$, can be computed in exactly the same fashion. Next we consider the sum modulo $q_2$, which is given by
$$
\mathcal T_2=\sum_{\alpha\bmod{q_2}}\;\mathop{\sideset{}{^\star}\sum\sideset{}{^\star}\sum}_{\beta,\gamma\bmod{q_2}}e_{q_2}(\bar{q_1}h\beta-\bar{q_1}n\bar{\beta}-\bar{\tilde q}_1h\gamma+\bar{\tilde q}_1n\bar{\gamma})S(\bar{q_1}\bar{\beta},\bar{q_1}\alpha;q_2)S(\bar{\tilde q}_1\bar{\gamma},\bar{\tilde q}_1\alpha;q_2)e_{q_2}\left(\bar{\tilde q}_1\bar{q_1}m\alpha\right).
$$
Opening the Kloosterman sums we execute the sum over $\alpha$ to get
$$
\mathcal T_2=q_2\sideset{}{^{\star\star}}\sum_{\delta\bmod{q_2}}\;\mathop{\sideset{}{^\star}\sum\sideset{}{^\star}\sum}_{\beta,\gamma\bmod{q_2}}e_{q_2}(\bar{q_1}h\beta-\bar{q_1}n\bar{\beta}-\bar{\tilde q}_1h\gamma+\bar{\tilde q}_1n\bar{\gamma}+\bar{q_1}\bar{\beta}\bar{\delta}-\bar{\tilde q}_1q_1\bar{\gamma}(\overline{\tilde q_1\delta+m})),
$$
where the double asterisk on the sum over $\delta$ indicates that $\delta$ and $\tilde q_1\delta+m$ are invertible modulo $q_2$. To get square-root cancellation in the remaining character sum we shall appeal to the work of Bombieri and Sperger \cite{BS} (in particular see Section IV. of \cite{BS}). 
 
Using the notation of \cite{BS}, we set 
$$
f(x,y,z)=a(x)+y+\frac{b(x)}{y}+z+\frac{c(x)}{z}
$$  
where
$$
a(x)=0,\;\;\;b(x)=\bar{q}_1^2h\left(\frac{1}{x}-n\right),\;\;\;\text{and}\;\;\;c(x)=\bar{\tilde q}_1^2h\left(\frac{q_1}{\tilde q_1x+m}-n\right).
$$ 
Let $\mathcal V$ be the quasi-projective variety in $\mathbb A_3(\mathbb F_{q_2})$ defined by $x\neq 0, -\bar{\tilde q}_1m$, $y\neq 0$, and $z\neq 0$. Then
$$
\mathcal T_2=q_2\sum_{(x,y,z)\in \mathcal V(\mathbb F_{q_2})}e_{q_2}(f(x,y,z)).
$$
From the main result of \cite{BS}, it follows that if $q_2\nmid q_1-\tilde q_1$ or $q_2\nmid m$, then $\mathcal T_2\ll q_2^{\frac{5}{2}}$. Otherwise using the Weil bound for Kloosterman sums we have $\mathcal T_2\ll q_2^3$. Putting everything together we conclude the following bound.
\begin{lemma}
\label{char-sum-1}
For $q_1\neq \tilde q_1$, the character sum $\mathcal T(n,m,h;q_1,\tilde q_1,q_2)$ vanishes unless $(m,q_1\tilde q_1)=1$, in which case we have
$$
\mathcal T(n,m,h;q_1,\tilde q_1,q_2)\ll q_1^{\frac{3}{2}}\tilde q_1^{\frac{3}{2}}q_2^{\frac{5}{2}}(m,q_2)^{\frac{1}{2}}.
$$
\end{lemma} 

If $q_1=\tilde q_1$, then the character sum $\mathcal T(n,m,h;q_1,\tilde q_1,q_2)$ splits as a product of two character sums. The one with modulus $q_2$, has already been analysed above and it satisfies the bound $\ll q_2^{\frac{5}{2}}(m,q_2)^{\frac{1}{2}}$. The other sum with modulus $q_1$ is given by
$$
\mathcal T_1^{q_1=\tilde q_1}=\sum_{\alpha\bmod{q_1}}\mathop{\sideset{}{^\star}\sum\sideset{}{^\star}\sum}_{\beta,\gamma\bmod{q_1}}e_{q_1}(\bar{q_2}h\beta-\bar{q_2}n\bar{\beta}-\bar{q_2}h\gamma+\bar{q_2}n\bar{\gamma})S(\bar{q_2}\bar{\beta},\bar{q_2}\alpha;q_1)S(\bar{q_2}\bar{\gamma},\bar{q_2}\alpha;q_1)e_{q_1^2}\left(\bar{q_2}m\alpha\right).
$$
As before we open the Kloosterman sums and execute the sum over $\alpha$. It follows that the sum vanishes unless $q_1|m$. So let $m=q_1m'$. Then we arrive at
$$
\mathcal T_1^{q_1=\tilde q_1}=q_1\sideset{}{^{\star\star}}\sum_{\delta\bmod{q_1}}\mathop{\sideset{}{^\star}\sum\sideset{}{^\star}\sum}_{\beta,\gamma\bmod{q_1}}e_{q_1}(\bar{q_2}h\beta-\bar{q_2}n\bar{\beta}-\bar{q_2}h\gamma+\bar{q_2}n\bar{\gamma}+\bar{q_2}\bar{\beta}\bar{\delta}-\bar{q_2}
\bar{\gamma}\overline{(\delta+m')}).
$$ 
The remaining sum is just like $\mathcal T_2$, with different parameters. 
\begin{lemma}
\label{char-sum-2}
The character sum $\mathcal T(n,m,h;q_1,q_1,q_2)$ vanishes unless $q_1|m$, in which case we have
$$
\mathcal T(n,q_1m',h;q_1,q_1,q_2)\ll q_1^{\frac{5}{2}}q_2^{\frac{5}{2}}\sqrt{(m',q_1q_2)}.
$$
\end{lemma}

\subsection{Final estimation and conclusion of the theorem}

It follows from Lemma \ref{char-sum-1}, that
\begin{align}
\label{bd-char-1}
\frac{M^2X^3}{Q_1Q^4}\mathop{\sum_{q_1\in\mathcal Q_1}\sum_{\tilde q_1\in\mathcal Q_1}}_{q_1\neq \tilde q_1}\sum_{|m|\ll \frac{Q_1Q}{M}X^{\varepsilon}}|\mathcal T(n,m,h;q_1,\tilde q_1,q_2)|&\ll \frac{M^2X^3Q_1^5Q_2^{\frac{5}{2}}}{Q_1Q^4}\sum_{1\leq |m|\ll \frac{Q_1Q}{M}X^{\varepsilon}}\sqrt{(m,q_2)}\\
\nonumber &\ll \frac{MX^3Q_1^5Q_2^{\frac{5}{2}}}{Q^3}X^{\varepsilon}.
\end{align}
Applying Lemma \ref{char-sum-2}, it follows that
\begin{align}
\label{bd-char-2}
&\frac{M^2X^3}{Q_1Q^4}\mathop{\sum_{q_1\in\mathcal Q_1}\sum_{\tilde q_1\in\mathcal Q_1}}_{q_1=\tilde q_1}\sum_{|m|\ll \frac{Q_1Q}{M}X^{\varepsilon}}|\mathcal T(n,m,h;q_1,\tilde q_1,q_2)|\\
\nonumber =&\frac{M^2X^3}{Q_1Q^4}\sum_{q_1\in\mathcal Q_1}\sum_{1\leq |m|\ll \frac{Q}{M}X^{\varepsilon}}|\mathcal T(n,q_1m,h;q_1,q_1,q_2)|+\frac{M^2X^3}{Q_1Q^4}\sum_{q_1\in\mathcal Q_1}|\mathcal T(n,0,h;q_1,q_1,q_2)|\\
\nonumber \ll &\frac{MX^3}{\sqrt{Q}}X^{\varepsilon}+\frac{M^2X^3}{Q}X^{\varepsilon}\ll \frac{M^2X^3}{Q}X^{\varepsilon}.
\end{align}
The above two bounds \eqref{bd-char-1}, \eqref{bd-char-2}, together with Lemma \ref{rem-sum} imply that
\begin{align*}
\tilde D_{h,\alpha,1}^\sharp(X,M;n,q_2)\ll \frac{MX^3Q_1^2}{\sqrt{Q_2}}X^{\varepsilon}+\frac{M^2X^3}{Q_1Q_2}X^{\varepsilon}.
\end{align*}
Plugging this estimate in \eqref{dh} we get the following:
\begin{lemma}
\label{last-lemma}
For $Q_1Q_2=Q=X^{\frac{1}{2}+\delta}$ and $M\ll X^{\frac{1}{2}+3\delta+\varepsilon}$, we have 
\begin{align*}
\tilde D_{h,\alpha,1}^\sharp(X,M)\ll \left(\frac{X^{\delta}}{Q_2^{\frac{1}{4}}}+\frac{X^{2\delta}}{Q_1}\right)X^{1+\varepsilon}.
\end{align*}
\end{lemma}
For given $\delta$ the optimal break up of $Q$ is given by $Q_1=X^{\frac{1}{10}+\delta}$ and $Q_2=X^{\frac{2}{5}}$. To obtain the optimal value for $\delta$, recall that we are assuming $\min\{Q_1,Q_2\}\gg X^{2\delta+\varepsilon}$. So we have the restriction $\delta<\frac{1}{10}$. Moreover comparing the above bound with the bound for the error term in \eqref{exp1}, we get that the optimal choice is given by
$$
\delta=\frac{1}{10}-\delta
$$ 
i.e. $\delta=\frac{1}{20}$. Our main theorem now follows from Lemma \ref{first-lemma}, Lemma \ref{mid-lemma}, Lemma \ref{last-but-one} and Lemma \ref{last-lemma}.


\end{document}